\documentstyle[12pt]{article}

\advance \topmargin by -\headheight
\advance \topmargin by -\headsep
     
\evensidemargin \oddsidemargin
\begin{document}
 
\topmargin -.6in

\def\y {\'{\i}}
\def\beq{\begin{equation}}
\def\eeq{\end{equation}}
\def\bray{\begin{eqnarray}}
\def\eray{\end{eqnarray}}

\begin{titlepage}
${}$
January 1998 

\hfill{q-alg/xxxxxxx}
 \vskip .6in
 
\begin{center}
{\large {\bf Hamiltonian Reduction and the Construction of q-Deformed
Extensions of the Virasoro Algebra}}
\end{center}

\normalsize
\vskip .4in
 
\begin{center}

{E. Batista \footnotemark \footnotetext{Supported by FAPESP} ,
J. F. Gomes\footnotemark \footnotetext{Work partially
 supported by CNPq} , 
I.J. Lautenschleguer \footnotemark \footnotetext {Supported by CNPq}} \\
\par \vskip .1in \noindent
Instituto de F\'{\i}sica Te\'{o}rica-UNESP\\
Rua Pamplona 145\\
01405-900 S\~{a}o Paulo, S.P.\\
 Brazil \\
\par \vskip .3in

\end{center}
 
\begin{center}
{\large {\bf ABSTRACT}}
\end{center}
In this paper we employ the construction of Dirac bracket for the remaining
current of $sl(2)_q$ deformed Kac-Moody algebra when constraints similar
to those connecting the $sl(2)$-WZW model and the Liouville theory are
imposed and show that it satisfy the q-Virasoro algebra proposed by Frenkel
and Reshetikhin. The crucial assumption considered in our calculation is
 the existence
of a classical Poisson bracket algebra induced, in a consistent manner by the
correspondence principle, mapping the quantum generators into commuting objects
of classical nature preserving their algebra.

\par \vskip .3in \noindent

\end{titlepage}


The Virasoro algebra and its extensions have been understood to provide the algebraic
structure underlying conformally invariant models which includes  string theory and 2D
statistical models on the lattice.  On the other hand, quantum groups also play an
important role in the integrability properties of those models (see for
instance \cite{saleurzuber}, \cite{Alekseev}, \cite{Goddard}, \cite{Faddeev}). 
 It thus, seems natural to connect these two important subjects by constructing
a q-deformed version of the Virasoro algebra and its extensions.  This may prove useful in
establishing a q-deformed string model in the line of refs. \cite{cgk},
\cite{devega}.

The construction of q-deformed Virasoro algebra have been proposed using both, bosons and
fermions \cite{cp}. However  a connection with the classical
canonical structure is still unclear.  Frenkel and Reshetikhin \cite{Frenkel2} have
proposed a q-Virasoro algebra based on the q-deformation of a Miura transformation
involving  classical Poisson brackets.
The  Hamiltonian reduction provides a systematic 
procedure in constructing extensions of the Virasoro algebra by adding to the
spin 2, generators of higher spin. A typical example of such procedure connects
Wess-Zumino-Witten (WZW) model to the 2D Toda field theories. 
The latter arises when a consistent
set of constraints are implemented to the Kac-Moody currents describing the
WZW model associated to a Lie group G \cite{O'Raifeartaigh1} or to an infinite dimensional
Kac-Moody group $\hat G $ \cite{aratyn}.

A redefinition of the canonical Poisson brackets into Dirac brackets is
required in order for the equations of motion of the reduced model to be
consistent with those obtained from the remaining current algebra. Under the
Dirac bracket the spin one generators corresponding to the remaining currents
become the $W_n$ generators of higher spin defined according to an improved
energy momentum tensor (see \cite{O'Raifeartaigh1} for a rewiew).

For the q-deformed Kac-Moody algebras, although a canonical structure is still
unknown, their algebra is well established \cite{Drinfeld},\cite{Jimbo} and can be
constructed in terms of non commuting objects (quantum fields) \cite
{Odake},\cite{Frenkel1}.

In this paper we employ the construction of Dirac bracket for the remaining
current of $sl(2)_q$ deformed Kac-Moody algebra when constraints similar
to those connecting the $sl(2)$ WZW model and the Liouville theory are
imposed. The crucial assumption considered in our calculation is the existence
of a classical Poisson bracket algebra induced, in a consistent manner by the
correspondence principle, mapping the quantum generators into commuting objects
of classical nature preserving their algebra.  We show that the remaining algebra coincide
with the q-Virasoro algebra proposed by Frenke and Reshetikhin \cite{Frenkel2}


For $q=1$, the classical $sl(2)$ Poisson bracket algebra derived from the WZW
model \cite{Witten} is given in terms of formal power series by
\bray
\{H(z),H(w)\}&=&-ik\sum_{n\in{\bf Z}}n\left(\frac{w}{z}\right)^n,
\label{cor1}\\
\{H(z),E^\pm(w)\}&=&\mp i\sqrt{2}E^\pm(z)\sum_{n\in{\bf Z}}\left(\frac{w}{z}\right)^n,
\label{cor2}\\
\{E^+(z),E^-(w)\}&=&-i\sqrt{2}H(z)\sum_{n\in{\bf Z}}\left(\frac{w}{z}\right)^n
-ik\sum_{n\in{\bf Z}}n\left(\frac{w}{z}\right)^n.
\label{cor3}
\eray
where $k$ characterizes the central term.
The corresponding conformal Toda model associated to $G=sl(2)$ (Liouville model) is
obtained by constraining \cite{O'Raifeartaigh1}
\beq
\chi_1=H(z)\approx 0\quad,\quad \chi_2=E^+(z)-1\approx 0.
\eeq

The Dirac bracket is defined by
\begin{eqnarray}
&\{A(z),B(w)\}_D&=\{A(z),B(w)\}_P  \nonumber \\
& &-\oint\oint\frac{du}{2\pi iu}\frac{dv}{2\pi iv}
\{A(z),\chi_i(u)\}\Delta^{-1}_{ij}(u,v)\{\chi_j(v),B(w)\}\nonumber\\
&&
\label{dirb}
\end{eqnarray}
where $\Delta^{-1}(x,y)$ is the inverse of the Dirac matrix $\Delta_{ij}(x,y)=
\{\chi_i(x),\chi_j(y)\}$ in the sense that
\beq
\oint\frac{du}{2\pi iu}\Delta_{ij}(z,u)\Delta^{-1}_{jk}(u,w)=\delta_{ik}
\sum_{n\in{\bf Z}}\left(\frac{z}{w}\right)^{n}=\delta_{ik}\delta\left(\frac{z}{w}\right).
\label{inver}
\eeq

Under the Dirac bracket, the remaining current $E^-(z)$ with $k=1$ 
leads to the Virasoro algebra
\beq
\{E^-(z),E^-(w)\}_D=-i(E^-(z)+E^-(w))\sum_{n\in{\bf Z}}n\left(\frac{w}{z}\right)^n
+\frac{i}{2}\sum_{n\in{\bf Z}}n^3\left(\frac{w}{z}\right)^n.
\eeq


We now consider the q-deformed Kac-Moody algebra for $sl(2)_q$ of level $k$ defined
by \cite{Jimbo},\cite{Drinfeld}
\bray
[H_n,H_m]&=&\frac{[2n][kn]}{2n}\delta_{m+n,0},\\
\,[H_0,H_m]&=&0,\\
\,[H_n,E^\pm_m]&=&\pm\sqrt{2}q^{\mp|n|\frac{k}{2}}\frac{[2n]}{2n}E^\pm_{m+n},\\
\,[H_0,E^\pm_m]&=&\pm\sqrt{2}E^\pm_m,\\
\,[E^+_n,E^-_m]&=&\frac{q^{\frac{k(n-m)}{2}}\Psi_{n+m}-q^{\frac{k(m-n)}{2}}\Phi_{n+m}}
{q-q^{-1}},\\
E^\pm_{n+1}E^\pm_m-q^{\pm 2}E^\pm_mE^\pm_{n+1}&=&q^{\pm 2}
 E^\pm_{n}E^\pm_{m+1}-
E^\pm_{m+1}E^\pm_{n}
\eray
where
\bray
\Psi(z)&=&q^{\sqrt{2}H_0}e^{\sqrt{2}(q-q^{-1})\sum_{n>0}H_nz^{-n}},
\label{mai}\\
\Phi(z)&=&q^{-\sqrt{2}H_0}e^{-\sqrt{2}(q-q^{-1})\sum_{n<0}H_nz^{-n}},
\label{men}
\eray
and $[x] = {{q^x - q^{-x} }\over {q - q^{-1}}}$, leading to the Operator Product relations
\bray
H(z)H(w)&=&\sum_{n>0}\frac{[2n][kn]}{2n}\left(\frac{w}{z}\right)^n,\\
H(z)E^\pm(w)&=&\pm\sqrt{2}\left(1+\sum_{n>0}\frac{[2n]}{2n}\left(\frac{wq^{\mp
{k\over 2}}}
{z}\right)^n\right)E^\pm(w),\\
E^+(z)E^-(w)&=&\frac{1}{w(q-q^{-1})}\left(\frac{\Psi(wq^{\frac{k}{2}})}{z-wq^k}-
\frac{\Phi(wq^{-\frac{k}{2}})}{z-wq^{-k}}\right),
\label{mame}\\
E^\pm(z)E^\pm(w)(z-wq^{\pm 2})&=&E^\pm(w)E^\pm(z)(zq^{\pm 2}-w).
\label{ncom}
\eray
for $|z| > |w|$ and we are considering $q$ to be a pure phase.
It is clear from (\ref{ncom}) that $E^\pm$ are not self commuting objects, however
this structure can be constructed using the Wakimoto construction \cite{Odake}. In
particular, for $k=1$, it can be constructed in terms of a single Fubini Veneziano
field \cite{Frenkel1} as follows
\beq
E^\pm(z)=:e^{\pm i\sqrt{2}Q^\pm(z)}:\quad,\quad H(z)=\sum_{n\in{\bf Z}}\alpha_nz^{-n},
\label{camp1}
\eeq
where
\beq
Q^\pm(z)=\tilde{q}-i\tilde{p}\ln{z}+i\sum_{n<0}\frac{\alpha_n}{[n]}(zq^{\mp
\frac{1}{2}})^{-n}+i\sum_{n>0}\frac{\alpha_n}{[n]}(zq^{\pm
\frac{1}{2}})^{-n},
\label{camp2}
\eeq
and
\beq
[\alpha_n,\alpha_m]=\frac{[2n][n]}{2n}\delta_{m+n,0}\quad;\quad[\tilde{q},\tilde{p}]
=i.
\label{camp3}
\eeq

We should point out that for $q=1$, the vertex operator construction (\ref{camp1})-(\ref
{camp3}) satisfy  (\ref{cor1})-(\ref{cor3}) with $k=1$. Our
Hamiltonian reduction procedure consist in
implementing the following constraints
\bray
\chi^q_1&=&\frac{\Psi(z)-\Phi(z)}{\sqrt{2}(q-q^{-1})}\approx 0,
\label{ain1}\\
\chi_2^q&=&E^+(z)\approx 1,
\label{ain2}
\eray
for $\Psi (z)$ and $\Phi(z)$ defined in (\ref{mai}) and (\ref{men}) respectively. Notice
that $\chi_1^q=H(z)+O(q-q^{-1})$, and reduce consistently to the known $q=1$ case.

For $q \neq 1$, the q-deformed Dirac matrix is constructed out of the following relations
obtained by direct calculation using the vertex operators (\ref{camp1})-(\ref{camp3})
\bray
\Psi(z)\Phi(w)&=&\frac{(z-wq^3)(z-wq^{-3})}{(z-wq)(z-wq^{-1})}\Phi(w)\Psi(z),
\label{ope1}\\
\Psi(z)E^\pm(w)&=&q^{\pm 2}\frac{(z-wq^{\mp\frac{5}{2}})}
{(z-wq^{\pm\frac{3}{2}})}E^\pm(w)\Psi(z),
\label{ope2}\\
E^\pm(z)\Phi(w)&=&q^{\pm 2}\frac{(z-wq^{\mp\frac{5}{2}})}
{(z-wq^{\pm\frac{3}{2}})}\Phi(w)E^\pm(z),
\label{ope3}
\eray
(for $|z|>|w|$) together with (\ref{mame}) and (\ref{ncom}) for $k=1$.

From equations (\ref{ope1})-(\ref{ope3}) we evaluate
\bray
\left[\frac{\Psi(z)-\Phi(z)}{\sqrt{2}(q-q^{-1})},\frac{\Psi(w)-\Phi(w)}{\sqrt{2}(q-q^{-1})}
\right]&=&\frac{[2]}{2}\Phi(w)\Psi(z)\sum_{n>0}\left(\frac{w}{z}\right)^n[n]\nonumber\\
&-&\frac{[2]}{2}\Phi(z)\Psi(w)\sum_{n>0}\left(\frac{z}{w}\right)^n[n],
 \label{eva1}
\eray
\bray
\left[\frac{\Psi(z)-\Phi(z)}{\sqrt{2}(q-q^{-1})},E^\pm(w)
\right]&=&\pm\frac{[2]}{\sqrt{2}}E^\pm(w)\Psi(z)\left(\sum_{n\geq 0}\left(\frac{wq^{\pm\frac{3}{2}}}
{z}\right)^n-\frac{q^{\mp 1}}{[2]}\right)\nonumber\\
&\pm&\frac{[2]}{\sqrt{2}}\Phi(z)E^\pm(w)\left(\sum_{n\geq 0}\left(\frac{zq^{\pm\frac{3}{2}}}
{w}\right)^n-\frac{q^{\mp 1}}{[2]}\right).
\label{eva2}
\eray
\newpage
\bray
[E^\pm(z),E^\mp(w)]&=&\pm\frac{1}{q-q^{-1}}\left(\Psi(wq^{\pm\frac{1}{2}})
\sum_{n\in{\bf Z}}\left(\frac{wq^{\pm 1}}{z}\right)^n\right.\nonumber\\
&-&\left.\Phi(wq^{\mp\frac{1}{2}})
\sum_{n\in{\bf Z}}\left(\frac{wq^{\mp 1}}{z}\right)^n\right)
\eray
\label{eva}
and
\bray
[E^\pm(z),E^\pm(w)]&=&
\pm\frac{1}{2}(q-q^{-1})E^\pm(w)E^\pm(z)\left([2]\sum_{n\geq 0}\left(\frac{wq^{\pm 2}}
{z}\right)^n-q^{\mp 1}\right)\nonumber \\
&\mp&\frac{1}{2}(q-q^{-1})E^\pm(z)E^\pm(w)\left([2]\sum_{n\geq 0}\left(\frac{zq^{\pm 2}}
{w}\right)^n-q^{\mp 1}\right).
\label{eva3}
\eray

Notice that the r.h.s. of (\ref{eva1}) e (\ref{eva2}) is normal ordered and all brackets display
explicit antisymmetry under $z\leftrightarrow w$.

Let us now discuss the classical counterpart of the quantum brackets (\ref{eva1})-(\ref{eva3}). The
usual canonical quantization procedure associates the classical Poisson bracket structure to quantum
commutators as
\beq
\{\;,\;\}\rightarrow -i[\;,\;].
\label{pres}
\eeq
The new feature compared with the $q=1$ case is the nonvanishing of eqn.
(\ref{eva3}).  Moreover, eqn. (\ref{eva3}) present a quadratic structure which
suggests an exponential realization (vertex operator for $k=1$ or the 
generalized Wakimoto construction for
generic $k$ (see \cite{Odake}))and the commutators are evaluated using the
 Baker-Haussdorff formula. The
latter has no classical analog but still expect a classical counterpart for
 the quantum algebra
(\ref{eva1})-(\ref{eva3}) to preserve their structure of algebraic nature.

We propose a classical Poisson bracket algebra by mapping quantum operators 
$\hat{A},\hat{B}$ into
classical objects $A,B$ such that
\beq
\{A(z),B(w)\}_{PB}\rightarrow \mp i[\hat{A}(z),\hat{B}(w)],
\label{princ}
\eeq
where the plus  sign is only taken for  $A=B = E^{\pm}$.  All other brackets
follow the usual correspondence principle (\ref{pres}).
Under this prescription and constraints (\ref{ain1}) and (\ref{ain2}), we construct
 the Dirac matrix
$\Delta_{ij}(z,w)=\{\chi_i(z),\chi_j(w)\}_{PB}$ to be
\bray
\Delta_{11}(z,w)&=&-i\frac{[2]}{2}\left(\sum_{n>0}\left(\frac{w}{z}\right)^n[n]-
\sum_{n>0}\left(\frac{z}{w}\right)^n[n]\right),
\label{elem1}\\
\Delta_{12}(z,w)&=&-i\frac{[2]}{\sqrt{2}}\left(\sum_{n\geq 0}\left(\frac{wq^{\frac{3}{2}}}{z}\right)^n+
\sum_{n\geq 0}\left(\frac{zq^{\frac{3}{2}}}{w}\right)^n\right)+\frac{2q^{-1}i}{\sqrt{2}},
\label{elem2}\\
\Delta_{22}(z,w)&=&i\frac{[2]}{2}(q-q^{-1})\left(\sum_{n>0}\left(\frac{wq^2}{z}\right)^n-
\sum_{n>0}\left(\frac{zq^2}{w}\right)^n\right),
\label{elem3}
\eray
and $\Delta_{21}(z,w)=-\Delta_{12}(w,z)$.

Its inverse is defined by equation (\ref{inver}) yielding
\bray
\Delta_{11}^{-1}(z,w)&=&\frac{-2i(q-q^{-1})}{[2]}\left(\sum_{n>0}\left(\frac{w}{z}\right)^n\frac{[n]}
{[2n]}-\sum_{n>0}\left(\frac{z}{w}\right)^n\frac{[n]}{[2n]}\right),\\
\Delta_{12}^{-1}(z,w)&=&\frac{-2i\sqrt{2}}{[2]}\left(\sum_{n>0}\left(\frac{wq^{-\frac{1}{2}}}{z}
\right)^n\frac{[n]}{[2n]}+\sum_{n>0}\left(\frac{zq^{-\frac{1}{2}}}{w}\right)^n
\frac{[n]}{[2n]}\right),\\
\Delta_{21}^{-1}(z,w)&=&\frac{2i\sqrt{2}}{[2]}\left(\sum_{n>0}\left(\frac{wq^{-\frac{1}{2}}}{z}
\right)^n\frac{[n]}{[2n]}+\sum_{n>0}\left(\frac{zq^{-\frac{1}{2}}}{w}\right)^n
\frac{[n]}{[2n]}\right),\\
\Delta_{22}^{-1}(z,w)&=&\frac{2i}{[2]}\left(\sum_{n>0}\left(\frac{wq^{-2}}{z}\right)^n\frac{[n]^2}{[2n]}-
\sum_{n>0}\left(\frac{zq^{-2}}{w}\right)^n\frac{[n]^2}{[2n]}\right),
\eray
and the Dirac bracket (\ref{dirb}) for the remaining current $E^-(z)$ can be evaluated
using the modified correspondence principle (\ref{princ}) in eqns. 
(\ref{eva1})-(\ref{eva3}) yielding, after
redefining $\tilde E^{-} = (q-q^{-1})^2\sqrt{\frac{[2]}{2}}
E^{-} +\frac{4}{\sqrt{2[2]}} $
\bray
\{\tilde E^-(z),\tilde E^-(w)\}_D&=&\frac{i[2]}{2}(q-q^{-1})^2
\tilde E^-(z)\tilde E^-(w)\sum_{n\in{\bf Z}}q^{-2|n|}
\frac{[n]^2}{[2n]}\left(\frac{z}{w}\right)^n\nonumber\\
&-&i(q-q^{-1})^2\sum_{n\in{\bf Z}}q^{-2|n|}[2n]\left(\frac{z}{w}\right)^n.
\label{qdirb}
\eray

The algebra given in (\ref{qdirb}) coincide, apart from the factor $q^{-2|n|}$ to the q-Virasoro algebra
proposed by Frenkel and Reshetikhin (see \cite{Frenkel2}) with $q=e^{ih}$. This undesirable
 factor may be absorbed by redefining the
classical brackets (\ref{elem1})-(\ref{elem3}) of the form
$$
\{A(z),B(w)\}=\sum_{n\in{\bf Z}}C_n\left(\frac{z}{w}\right)^n
$$
into
$$
\{A(z),B(w)\}=\sum_{n\in{\bf Z}}C_n\left(\frac{z}{w}\right)^nq^{2|n|}.
$$

Under this modification of the correspondence principle (\ref{princ}) the
 Dirac bracket for the remaining
current $\tilde E^-(z)$ coincide precisely with the algebra given in
\cite{Frenkel2}, namely,
\bray
\{\tilde E^-(z),\tilde E^-(w)\}_D&=&\frac{i[2]}{2}(q-q^{-1})^2
\tilde E^-(z)\tilde E^-(w)\sum_{n\in{\bf Z}}
\frac{[n]^2}{[2n]}\left(\frac{z}{w}\right)^n\nonumber\\
&-&i(q-q^{-1})^2\sum_{n\in{\bf Z}}[2n]\left(\frac{z}{w}\right)^n.
\label{qvir}
\eray


The differing factor $q^{2|n|}$ is viewed of quantum origin. It is known, for instance, in quantizing
the $SU(2)$ WZW model, that the coupling constant of the diagonal fields is shifted by a factor 2 (Coxeter
number of $SU(2)$ (see \cite{Alekseev},\cite{Goddard}). In geral we expect the 
quantum correction
associated to an q-deformed Kac-Moody algebra $\hat{g}$ to be given as
$$
\{A(z),B(w)\}=\sum_{n\in{\bf Z}}C_n\left(\frac{z}{w}\right)^nq^{h|n|}
$$ where $h$ is the Coxeter element of $g$.  For the general q-deformed Kac-Moody
algebra $\hat g $, if we follow the usual constraints connecting the $g$ invariant
WZW and the conformal Toda models \cite{O'Raifeartaigh1} we obtain the q-deformed
$W_n$-algebra by adding to the q-Virasoro (\ref{qvir}), generators of higher spin. 
The 
 construction of  a classical action invariant under
transformations generated by operators satisfying the proposed classical Poisson
algebra is also an interesting problem that are under investigation and shall be
reported in a future publication. 

{\large {\bf Aknowledgements}:}
 We thank Prof. A. H. Zimerman for many helpful discussions.

\end{document}